\input amstex
\voffset=-.5in

\NoBlackBoxes
\magnification=\magstep1
\documentstyle{amsppt}

\rightheadtext {Lefschetz number}

\topmatter

\title The Lefschetz-Hopf theorem and axioms for the Lefschetz number
\endtitle

\author Martin Arkowitz and Robert F. Brown \endauthor

\affil Dartmouth College, Hanover and University of California, Los 
Angeles \endaffil


\address
Hanover,  NH 03755-1890, USA
\endaddress
\email Martin.A.Arkowitz\@Dartmouth.edu
\endemail

\address
Los Angeles, CA 90095-1555, USA
\endaddress
\email
rfb\@math.ucla.edu,
\endemail

\abstract
The reduced Lefschetz number, that is, $L(\cdot) - 1$ where $L(\cdot)$ 
denotes the Lefschetz number, is proved to be the unique integer-valued 
function $\lambda$ on selfmaps of compact polyhedra which is constant 
on homotopy classes such that (1) $\lambda (fg) = \lambda (gf)$, for $f 
\colon X \to Y$ and $g \colon Y \to X$; (2) if $(f_1, f_2, f_3)$ is a 
map of a cofiber sequence into itself, then $\lambda (f_1) = \lambda 
(f_1) + \lambda (f_3)$; (3) $\lambda (f) = - (\deg(p_1 f e_1) + \dots + 
\deg(p_k f e_k))$, where $f$ is a selfmap of a wedge of $k$ circles, 
$e_r$ is the inclusion of a circle into the $r$th summand and $p_r$ is 
the projection onto the $r$th summand.  If $f \colon X \to X$ is a 
selfmap of a polyhedron and $I(f)$ is the fixed point index of $f$ on 
all of $X$, then we show that $I(\cdot) - 1$ satisfies the above 
axioms.  This gives a new proof of the Normalization Theorem: If $f 
\colon X \to X$ is a selfmap of a polyhedron, then $I(f)$ equals the 
Lefschetz number $L(f)$ of $f$.  This result is equivalent to the 
Lefschetz-Hopf Theorem: If $f \colon X \to X$ is a selfmap of a finite 
simplicial complex with a finite number of fixed points, each lying in 
a maximal simplex, then the Lefschetz number of $f$ is the sum of the 
indices of all the fixed points of $f$.
\endabstract

\endtopmatter

\document

\head 1. Introduction. \endhead

Let $X$ be a finite polyhedron and denote by $\widetilde H_*(X)$ its
reduced homology with rational coefficients.  Then the {\it reduced
Euler
characteristic} of $X$, denoted by $\tilde \chi (X)$, is defined by
$$
\tilde \chi (X) = \sum_j (-1)^j \dim \widetilde H_j(X).
$$
Clearly, $\tilde \chi (X)$ is just the Euler characteristic minus one.
In 1962, Watts \cite{13} characterized the reduced Euler
characteristic as follows:  Let $\epsilon$ be a function from the set
of finite polyhedra with base points to the integers such that (i)
$\epsilon(S^0) = 1$, where
$S^0$ is the $0$-sphere, and (ii) $\epsilon(X) = \epsilon(A) +
\epsilon(X/A)$, where $A$ a subpolyhedron of $X$.  Then $\epsilon(X)
= \tilde \chi(X)$.

Let $\Cal C$ be the collection of spaces $X$ of the homotopy type of a
finite, connected CW-complex. If $X \in \Cal C$, we do not assume that
$X$ has a base point except
when $X$ is a sphere or a wedge of spheres.  It is not
assumed that maps between spaces with base points are based.
A map $f \colon X \to X$, where $X \in
\Cal C$, induces trivial homomorphisms $f_j \colon
H_j(X) \to H_j(X)$ of rational homology vector spaces for
all $j > \dim X$.
The {\it Lefschetz number}
$L(f)$ of $f$ is defined by
$$
L(f) = \sum_j (-1)^j Tr\, f_j,
$$
where $Tr$ denotes the trace.
The reduced Lefschetz number $\widetilde L$ is given by $\widetilde
L(f) = L(f) - 1$ or, equivalently, by considering the rational,
reduced homology homomorphism induced by $f$.

Since $\widetilde L(id) = \tilde \chi(X)$, where $id \colon X \to X$
is the identity map, Watts's Theorem suggests an axiomatization for
the reduced Lefschetz number which we state below as Theorem 1.1.

For $k \ge 1$, denote by  $\bigvee^k S^n$ the
wedge of $k$ copies of the $n$-sphere $S^n, n \ge 1$.  If we write
$\bigvee^k S^n$
as $S^n_1 \vee S^n_2 \vee \cdots \vee S^n_k$, where $S^n_j = S^n$,
then we have inclusions $e_j \colon
S^n_j \to \bigvee^k S^n$
into the $j$-th summand and
projections
$p_j \colon \bigvee^k S^n \to S^n_j$
onto the $j$-th summand, for $j = 1,
\dots, k$.  If $f \colon \bigvee^k S^n \to \bigvee^k S^n$ is a
map, then $f_j \colon S^n_j \to S^n_j$ denotes the composition
$p_j f e_j$.  The
degree of a map $f \colon S^n \to S^n$ is denoted by $\deg(f)$.

We characterize the reduced Lefschetz number as follows.

\proclaim{Theorem 1.1} The reduced Lefschetz number $\widetilde L$ is
the unique
function $\lambda$  from the set of self-maps of spaces in $\Cal C$
to the integers that satisfies the following conditions:

1. (Homotopy Axiom) If $f,g \colon X \to X$ are homotopic maps, then
$\lambda(f) = \lambda(g)$.

2. (Cofibration Axiom) If
$A$ is a subpolyhedron of $X$,
$A \to X \to X/A$ is the resulting cofiber sequence
and there exists a commutative diagram
$$
\CD
A @>>> X @>>> X/A\\
@V f' VV @V f VV @V \bar f VV\\
A @>>> X @>>> X/A,
\endCD
$$
then $\lambda(f) = \lambda(f') + \lambda(\bar f)$.

3. (Commutativity Axiom) If $f \colon
X \to Y$ and $g \colon Y \to X$ are maps, then $\lambda(gf) =
\lambda(fg)$.

4. (Wedge of Circles Axiom) If $f \colon
\bigvee^k S^1 \to \bigvee^k S^1$
is a map, $k \ge 1$, then
$$
\lambda(f) = -(\deg(f_1) + \cdots + \deg(f_k)),
$$
where $f_j = p_j f e_j$.
\endproclaim

In an unpublished dissertation \cite{10}, Hoang extended Watts's
axioms to characterize the reduced Lefschetz number for
basepoint-preserving self-maps of finite polyhedra.
His list of axioms is different from, but similar to, those in
Theorem 1.1.

One of the classical results of fixed point theory is

\proclaim{Theorem 1.2 (Lefschetz-Hopf)} If $f \colon X \to X$ is
a map of a finite polyhedron with a finite set of fixed points, each
of which lies in a maximal simplex of $X$, then $L(f)$ is the sum of
the indices of all the fixed points of $f$. \endproclaim

The history of this result is described in \cite{3}, see also \cite{8, 
p.\,458}.   A proof that depends on a delicate argument due to Dold 
\cite{5} can be found in
\cite{2} and, in a more condensed form, in \cite{4}.  In an
appendix to his dissertation \cite{12}, D. McCord outlined a possibly
more direct argument, but no details were published.  The book of 
Granas and Dugundji \cite{8, pp.\, 441 - 450} presents an argument 
based on classical techniques of Hopf \cite{11}.  We use the
characterization of the reduced Lefschetz number in Theorem 1.1 to
prove the Lefschetz-Hopf theorem in a quite natural manner by showing
that the fixed point index satisfies the axioms of Theorem 1.1.
That is, we prove

\proclaim{Theorem 1.3 (Normalization Property)} If $f \colon X \to X$
is any map of a finite polyhedron, then $L(f) = i(X, f, X)$, the fixed
point index of $f$ on all of $X$. \endproclaim

The Lefschetz-Hopf Theorem follows from the Normalization Property by
the Additivity Property of the fixed point index.  In fact these
two statements are equivalent.  The Hopf Construction \cite{2, p.\,117}
implies that a map $f$ from a finite polyhedron to itself is homotopic
to a map that satisfies the hypotheses of the
Lefschetz-Hopf theorem.  Thus the Homotopy and Additivity Properties of
the fixed point index imply that the Normalization Property follows
from the Lefschetz-Hopf Theorem.

\head 2. Lefschetz numbers and exact sequences. \endhead

In this section, all vector spaces are over a fixed field $F$, which
will not be
mentioned, and are finite dimensional.  A graded vector space $V =
\{V_n\}$ will always have the following properties: (1) each $V_n$ is
finite dimensional and (2) $V_n = 0$ for $n < 0$ and for $n > N$, for
some non-negative integer $N$.  A {\it
map} $f \colon V \to W$ of graded vector spaces $V = \{V_n\}$ and $W
= \{W_n\}$ is a sequence of linear transformations $f_n \colon V_n
\to W_n$.  For a map $f \colon V \to V$, the {\it Lefschetz number}
is defined by
$$
L(f) = \sum_n (-1)^n Tr\, f_n.
$$

The proof of the following lemma is straightforward, and hence omitted.

\proclaim{Lemma 2.1} Given a map of short exact sequences of vector
spaces
$$
\CD
0 @>>> U @>>> V @>>> W @>>> 0\\
@. @V f VV @V g VV @V h VV @.\\
0 @>>> U @>>> V @>>> W @>>> 0,
\endCD
$$
then $Tr\,g = Tr\,f + Tr\,h$.
\qed \endproclaim

\proclaim{Theorem 2.2} Let $A, B$ and $C$ be graded vector spaces with
maps $\alpha \colon A \to B, \beta \colon B \to C$ and selfmaps
$f \colon A \to A, g \colon B \to B$ and $h \colon C \to C$.  If for
every $n$, there is a linear transformation $\partial_n \colon C_n \to
A_{n-1}$ such that the following diagram is commutative and has exact
rows:
$$
\CD
0 @>>> A_N @>\alpha_N>> B_N @>\beta_N>> C_N @>\partial_N>> A_{N-1}
@>\alpha_{N-1}>> \cdots\\
@. @V f_N VV @V g_N VV @V h_N VV @V f_{N-1} VV\\
0 @>>> A_N @>\alpha_N>> B_N @>\beta_N>> C_N @>\partial_N>> A_{N-1}
@>\alpha_{N-1}>> \cdots
\endCD
$$
$$
\CD
\cdots @>\partial_1>> A_0 @>\alpha_0>> B_0 @>\beta_0>> C_0 @>>> 0\\
@. @V f_0 VV @V g_0 VV @V h_0 VV @.\\
\cdots @>\partial_1>> A_0 @>\alpha_0>> B_0 @>\beta_0>> C_0 @>>> 0,
\endCD
$$
then
$$
L(g) = L(f) + L(h).
$$
\endproclaim

\demo{Proof} Let $Im$ denote the image of a linear transformation and
consider the commutative diagram
$$
\CD
0 @>>> Im\, \beta_n @>>> C_n @>>> Im\, \partial_n @>>> 0\\
@. @V h_n|Im\, \beta_n VV @V h_n VV @V f_{n-1}|Im\, \partial_n VV @.\\
0 @>>> Im\, \beta_n @>>> C_n @>>> Im\, \partial_n @>>> 0.
\endCD
$$
By Lemma 2.1, $Tr (h_n) = Tr (h_n|Im\, \beta_n) + Tr (f_{n-1}|Im\,
\partial_n)$.
Similarly, the commutative diagram
$$
\CD
0 @>>> Im\, \partial_n @>>> A_{n-1} @>>> Im\, \alpha_{n-1} @>>> 0\\
@. @V f_{n-1}|Im\, \partial_n VV @V f_{n-1} VV @V g_{n-1}|Im\,
\alpha_{n-1} VV @.\\
0 @>>> Im\, \partial_n @>>> A_{n-1} @>>> Im\, \alpha_{n-1} @>>> 0
\endCD
$$
yields $Tr (f_{n-1}|Im\, \partial_n) = Tr (f_{n-1}) - Tr (g_{n-1}|Im\,
\alpha_{n-1})$.  Therefore
$$
Tr (h_n) = Tr (h_n|Im\, \beta_n) + Tr (f_{n-1}) - Tr (g_{n-1}|Im\,
\alpha_{n-1}).
$$
Now consider
$$
\CD
0 @>>> Im\, \alpha_{n-1} @>>> B_{n-1} @>>> Im\, \beta_{n-1} @>>> 0\\
@. @V g_{n-1}|Im\, \alpha_{n-1} VV @V g_{n-1} VV @V h_{n-1}|Im\,
\beta_{n-1} VV @.\\
0 @>>> Im\, \alpha_{n-1} @>>> B_{n-1} @>>> Im\, \beta_{n-1} @>>> 0,
\endCD
$$
so $Tr (g_{n-1}|Im\, \alpha_{n-1}) = Tr (g_{n-1}) - Tr (h_{n-1}|Im\,
\beta_{n-1})$.   Putting this all together, we obtain
$$
Tr (h_n) = Tr (h_n|Im\, \beta_n) + Tr (f_{n-1}) - Tr (g_{n-1}) +
Tr (h_{n-1}|Im\, \beta_{n-1}).
$$
We next look at the left end of the original diagram and get
$$
0 = Tr (h_{N+1}) = Tr (f_N) - Tr (g_N) + Tr (h_N|Im\, \beta_N)
$$
and at the right end which gives
$$
Tr (h_1) = Tr (h_1|Im\, \beta_1) + Tr (f_0) - Tr (g_0) + Tr (h_0).
$$
A simple calculation now yields
$$
\align
\sum^N_{n = 0} (-1)^nTr(h_n) =&\, \sum^{N+1}_{n = 0} =
(-1)^n(Tr(h_n|Im\,
\beta_n) + Tr(f_{n-1}) - Tr(g_{n-1})\\
& \qquad \qquad + Tr(h_{n-1}|Im\, \beta_{n-1}))\\
=&\, - \sum^N_{n = 0} (-1)^nTr(f_n) + \sum^N_{n = 0} =
(-1)^nTr(g_n).
\endalign
$$
Therefore $L(h) = -L(f) + L(g)$.
\qed \enddemo

We next give some simple consequences of Theorem 2.2.

\medskip

If $f \colon (X, A) \to (X, A)$ is a selfmap of a pair, where $X, A \in
\Cal C$, then $f$ determines $f_X \colon X \to X$ and $f_A \colon A
\to A$.  The map $f$ induces homomorphisms $f_j \colon H_j(X, A) \to
H_j(X, A)$ of relative homology with coefficients in $F$.   The {\it
relative
Lefschetz number} $L(f; X , A)$ is defined by
$$
L(f; X, A) = \sum_j (-1)^j Tr f_j.
$$

Applying Theorem 2.2 to the homology exact sequence of the pair $(X,
A)$, we obtain

\proclaim{Corollary 2.3} If $f \colon (X, A) \to (X, A)$ is a map of
pairs, where $X, A \in \Cal C$, then
$$
L(f; X, A) = L(f_X) - L(f_A).
$$
\endproclaim

This result was obtained by Bowszyc \cite{1}.

\proclaim{Corollary 2.4} Suppose $X = P \cup Q$ where $X,
P, Q \in \Cal C$ and $(X; P, Q)$ is an proper triad \cite{6,
p.\,34}.
If $f \colon X \to X$
is a map such that $f(P) \subseteq P$
and $f(Q) \subseteq Q$
then, for $f_P, f_Q$ and
$f_{P \cap Q}$ the restrictions of $f$ to
$P, Q$ and $P \cap Q$ respectively, we
have
$$
L(f) = L(f_P) + L(f_Q) - L(f_{P \cap Q}).
$$
\endproclaim

\demo{Proof} The map $f$ and its restrictions induce a map of the
Mayer-Vietoris homology sequence \cite{6, p.\,39} to itself so the
result
follows from Theorem
2.2. \qed \enddemo

A similar result was obtained by Ferrario \cite{7, Theorem\,3.2.1}.

\medskip

Our final consequence of Theorem 2.2 will be used in the
characterization of the reduced Lefschetz number.

\proclaim{Corollary 2.5} If
$A$ is a subpolyhedron of $X$,
$A \to X \to X/A$ is the resulting cofiber sequence
of spaces in $\Cal C$ and there exists a commutative diagram
$$
\CD
A @>>> X @>>> X/A\\
@V f' VV @V f VV @V \bar f VV\\
A @>>> X @>>> X/A,
\endCD
$$
then
$$
L(f) = L(f') + L(\bar f) - 1.
$$
\endproclaim

\demo{Proof} We apply Theorem 2.2 to the homology cofiber sequence.
The `minus one' on the right hand side arises because that sequence
ends with
$$
\to H_0(A) \to H_0(X) \to \tilde H_0(X/A) \to 0. \quad \qed
$$
\enddemo

\head 3. Characterization of the Lefschetz number. \endhead

Throughout this section, all spaces are assumed to lie in $\Cal C$.

We let $\lambda$ be a function from the set of self-maps of spaces in
$\Cal C$ to the integers that satisfies the Homotopy Axiom,
Cofibration Axiom, Commutativity Axiom and Wedge of Circles Axiom of
Theorem 1.1 as stated in the Introduction.

We draw a few simple consequences of these axioms.  From the
Commutativity Axiom, we obtain

\proclaim{Lemma 3.1} If $f \colon X \to X$ is a map and $h \colon X
\to Y$ is a homotopy equivalence with homotopy inverse $k \colon Y
\to X$, then $\lambda(f) = \lambda(hfk)$. \qed \endproclaim

\proclaim{Lemma 3.2} If $f \colon X \to X$ is homotopic to a constant
map, then $\lambda(f) = 0$. \endproclaim

\demo{Proof} Let $*$ be a one-point space and $* \colon * \to *$ the
unique map.  From the map of cofiber sequences
$$
\CD
* @>>> * @>>> *\\
@V * VV @V * VV @V * VV\\
* @>>> * @>>> *
\endCD
$$
and the Cofibration Axiom, we have $\lambda(*) = \lambda(*) +
\lambda(*)$, and therefore $\lambda(*) = 0$.  Write any constant map
$c \colon X \to X$ as $c(x) = *$ for some $* \in X$,
let $e \colon * \to
X$ be inclusion
and $p \colon X \to *$ projection.  Then $c = ep$ and $pe = *$, and
so
$\lambda(c) = 0$ by the Commutativity Axiom.  The lemma follows from
the Homotopy Axiom. \qed \enddemo

If $X$ is a based space with base point $*$, i.e., a sphere or wedge of
spheres, then the cone and suspension of $X$ are defined by $CX = X
\times I/(X \times 1 \cup * \times I)$ and $\Sigma X = CX/(X \times 0)$,
respectively.

\proclaim{Lemma 3.3} If $X$ is a based space, $f \colon X \to X$ is a
based map and $\Sigma f \colon \Sigma X
\to \Sigma X$ is the suspension of $f$, then $\lambda(\Sigma f) =
-\lambda(f)$. \endproclaim

\demo{Proof} Consider the maps of cofiber sequences
$$
\CD
X @>>> CX @>>> \Sigma X\\
@V f VV @V Cf VV @V \Sigma f VV\\
X @>>> CX @>>> \Sigma X.
\endCD
$$
Since $CX$ is contractible, $Cf$ is homotopic to a constant
map.  Therefore, by Lemma 3.2 and the Cofibration Axiom,
$$
0 = \lambda(Cf) = \lambda(\Sigma f) + \lambda(f). \quad \qed
$$
\enddemo

\proclaim{Lemma 3.4} For any $k \ge 1$ and $n \ge 1$, if $f \colon
\bigvee^k S^n \to \bigvee^k S^n$ is a map, then
$$
\lambda(f) = (-1)^n(\deg(f_1) + \cdots + \deg(f_k)),
$$
where $e_r \colon S^n
\to \bigvee^k S^n$ and $p_r \colon \bigvee^k S^n \to S^n$ for $r = 1,
\dots, k$ are the inclusions and projections, respectively, and
$f_r = p_r f e_r$.
\endproclaim

\demo{Proof} The proof is by induction on the dimension $n$ of the
spheres.  The case $n = 1$ is the Wedge of Circles Axiom.  If $n \ge
2$, then the map $f \colon \bigvee^k S^n
\to \bigvee^k S^n$ is homotopic to a based map $f' \colon \bigvee^k S^n
\to \bigvee ^k S^n$.  Then $f'$ is homotopic to
$\Sigma g$, for some map $g
\colon \bigvee^k S^{n-1} \to \bigvee^k S^{n-1}$.  Note that if
$g_j \colon S^{n - 1}_j \to S^{n - 1}_j$, then $\Sigma g_j$ is
homotopic to $f_j \colon S^n_j \to S^n_j$.  Therefore by Lemma 3.3
and the induction hypothesis,
$$
\align
\lambda(f) = \lambda(f') = -\lambda(g) &=
-(-1)^{n-1}((\deg(g_1) +
\cdots +
\deg(g_k))\\
&= (-1)^n(\deg(f_1) + \cdots + \deg(f_k)).  \qed
\endalign
$$
\enddemo

\proclaim{Proof of Theorem 1.1}\endproclaim
\noindent
Since $\tilde L(f) = L(f) - 1$, Corollary 2.5 implies
that $\tilde L$ satisfies the Cofibration Axiom.
We next show that $\tilde L$ satisfies the Wedge of Circles Axiom.
There is an isomorphism $\theta \colon \bigoplus^k H_1(S^1) \to
H_1(\bigvee^k S^1)$ defined by $\theta(x_1, \dots , x_k) =
e_{1*}(x_1)\,
+\, \cdots \,+\, e_{k*}(x_k)$, where $x_i \in H_1(S^1)$.  The inverse
$\theta^{-1} \colon H_1(\bigvee^k S^1) \to \bigoplus^k H_1(S^1)$ is
given by $\theta^{-1}(y) = (p_{1*}(y), \dots , p_{k*}(y))$.
If $u \in H_1(S^1)$ is a generator, then a basis
for $H_1(\bigvee^k S^1)$ is $e_{1*}(u), \dots , e_{k*}(u)$.  By
calculating the trace of $f_* \colon H_1(\bigvee^k S^1) \to
H_1(\bigvee^k S^1)$ with respect to this basis, we obtain $\tilde L(f)
= -(\deg(f_1) + \cdots + \deg(f_k))$.
The remaining axioms
are obviously satisfied by $\tilde L$.  Thus $\tilde L$ satisfies the
axioms of Theorem 1.1.

Now suppose $\lambda$ is a
function from the self-maps of spaces in $\Cal C$ to the integers
that satisfies the axioms.  We regard $X$ as a
connected, finite
CW-complex and proceed by induction on the dimension of $X$.  If
$X$ is $1$-dimensional, then it is the homotopy type of a wedge of
circles.  By Lemma 3.1, we can regard $f$ as a self-map of $\bigvee^k
S^1$, and so the Wedge of Circles Axiom gives
$$
\lambda(f) = -(\deg(f_1) + \cdots + \deg(f_k)) = \tilde
L(f).
$$
Now suppose that $X$ is $n$-dimensional and let $X^{n-1}$ denote
the $(n - 1)$-skeleton of $X$.  Then $f$ is homotopic to a cellular
map $g \colon X \to X$ by the Cellular Approximation Theorem
\cite{9, Theorem\,4.8, p.\,349}.  Thus $g(X^{n-1}) \subseteq
X^{n-1}$, and so
we have a commutative diagram
$$
\CD
X^{n-1} @>>> X @>>> X/ X^{n-1} = \bigvee^k S^n\\
@V g' VV @V g VV @V \bar g VV\\
X^{n-1} @>>> X @>>> X/ X^{n-1} = \bigvee^k S^n.
\endCD
$$
Then, by the Cofibration Axiom, $\lambda(g) = \lambda(g') +
\lambda(\bar g)$.  Lemma 3.4 implies that $\lambda(\bar g) = \tilde
L(\bar g)$ so, applying the induction hypothesis to $g'$, we have
$\lambda(g) = \tilde L(g') + \tilde L(\bar g)$.  Since we have seen
that the reduced Lefschetz number satisfies the Cofibration Axiom, we
conclude that $\lambda(g) = \tilde L(g)$.  By the Homotopy Axiom,
$\lambda(f) = \tilde L(f)$.  \qed

\head 4. The Normalization Property. \endhead

Let $X$ be a finite polyhedron and $f \colon X \to X$ a map.  Denote
by $I(f)$ the fixed point index of $f$ on all of $X$, that is, $I(f)
= i(X, f, X)$ in the notation of \cite{2} and let $\tilde I(f) =
I(f) - 1$.

In this section we prove Theorem 1.3 by showing that, with rational
coefficients, $I(f) = L(f)$.

\proclaim{Proof of Theorem 1.3} \endproclaim

\noindent
We will prove that $\tilde I$ satisfies the axioms and
therefore, by Theorem 1.1, $\tilde I(f) = \tilde L(f)$.  The
Homotopy and Commutativity Axioms are well-known properties of the
fixed point index (see \cite{2, pp.\,59 and 62}).

To show that $\tilde I$ satisfies the Cofibration Axiom, it suffices
to consider $A$ a
subpolyhedron of $X$ and $f(A) \subseteq A$.  Let $f' \colon
A \to A$ denote the restriction of $f$ and $\bar f \colon X/A \to
X/A$ the map induced on quotient spaces.  Let $r \colon U \to A$
be a deformation retraction of a neighborhood of $A$ in $X$ onto $A$
and let $L$ be a
subpolyhedron of a barycentric subdivision of $X$ such that $A
\subseteq int \, L \subseteq L \subseteq U$.  By the Homotopy Extension
Theorem there is a homotopy $H \colon X \times I \to X$ such that
$H(x, 0) = f(x)$ for all $x \in X, H(a, t) = f(a)$ for all $a \in A$
and $H(x, 1) = fr(x)$ for all $x \in L$.  If we set $g(x) = H(x, 1)$
then, since there are no fixed points of $g$ on $L - A$, the
Additivity Property implies that
$$
I(g) = i(X, g, int \, L) + i(X, g, X - L).\leqno(4.1)
$$

We discuss each summand of (4.1) separately.  We begin with $i(X, g,
int \, L)$.
Since $g(L) \subseteq A \subseteq L$, it follows from the definition of
the index (\cite{2, p.\,56}) that $i(X, g, int \, L) = i(L, g, int \,
L)$.
Moreover, $ i(L, g, int \, L) = i(L, g, L)$ since there are no fixed
points on $L - int \, L$ (the Excision
Property of the index).  Let $e \colon A \to L$ be
inclusion then, by the Commutativity Property \cite{2, p.\,62} we have
$$
i(L, g, L) = i(L, eg, L) = i(A, ge, A) = I(f')
$$
because $f(a) = g(a)$ for all $a \in A$.

Next we consider the summand $i(X, g, X - L)$ of (4.1).
Let $\pi \colon X \to X/A$ be the quotient map, set $\pi(A) = *$ and
note that $\pi^{-1}(*) = A$.  If $\bar g \colon X/A \to X/A$ is
induced by $g$, the restriction of $\bar g$ to the
neighborhood $\pi(int \, L)$ of $*$ in $X/A$ is constant, so $i(X/A,
\bar g, \pi(int \, L)) = 1$.  If we denote the set of fixed points of
$\bar g$ with $*$ deleted by $Fix_* \bar g$, then $Fix_* \bar g$
is in the open subset $X/A - \pi(L)$ of $X/A$.  Let $W$ be an open
subset of $X/A$ such that $Fix_* \bar g \subseteq W \subseteq X/A -
\pi(L)$ with the property $\bar g(W) \cap \pi(L) = \emptyset$.  By
the Additivity Property we have
$$
I(\bar g) = i(X/A, \bar g, \pi(int \, L)) + i(X/A, \bar g, W) =
1 + i(X/A, \bar g, W).
$$
Now, identifying $X - L$ with the corresponding subset $\pi(X - L)$ of
$X/A$ and
identifying the restrictions of $\bar g$ and $g$ to those subsets, we
have $ i(X/A, \bar g, W) = i(X, g, \pi^{-1}(W))$.  The
Excision Property of the index implies that  $i(X, g, \pi^{-1}(W)) =
i(X, g, X - L)$.
Thus we have determined the second summand of (4.1): $i(X, g, X - L) =
I(\bar g) - 1$.

Therefore from (4.1) we obtain $I(g) = I(f') + I(\bar g) - 1$.  The
Homotopy
Property then tells us that
$$
I(f) = I(f') + I(\bar f) - 1
$$
since $f$ is homotopic to $g$ and $\bar f$ is homotopic to $\bar g$.
We conclude
that $\tilde I$ satisfies the Cofibration Axiom.

It remains to verify the Wedge of Circles Axiom.  Let $X =
\bigvee^k S^1 = S^1_1 \vee \cdots \vee S^1_k$
be a wedge of
circles with basepoint $*$ and $f \colon X \to X$ a map.
We first verify the axiom in the case $k = 1$.  We have $f \colon
S^1 \to S^1$ and we denote its degree by $deg(f) = d$.  We regard $S^1
\subseteq \Bbb C$, the complex numbers.  Then
$f$ is homotopic to $g_d$, where $g_d(z) = z^d$ has $|d - 1|$
fixed points for $d \ne 1$.  The fixed point index of $g_d$ in a
neighborhood of a
fixed point that contains no other fixed point of $g_d$ is $-1$ if $d
\ge 2$ and is $1$ if $d \le 0$.  Since $g_1$ is
homotopic to a map without fixed points, we see that $I(g_d) = -d +
1$ for all integers $d$.  We have shown that $I(f) = - deg(f) + 1$.

Now suppose $k \ge 2$.  If $f(*) = *$ then, by the Homotopy
Extension Theorem, $f$ is homotopic to a map which does not fix
$*$.  Thus we may assume, without loss of generality, that $f(*)
\in S^1_1 - \{*\}$.  Let $V$ be a neighborhood of $f(*)$ in
$S^1_1 -  \{*\}$ such that there exists a neighborhood
$U$ of $*$ in $X$ disjoint from $V$ with $f (\bar U) \subseteq
V$.  Since $\bar U$ contains no fixed point of $f$ and the open
subsets $S^1_j - \bar U$ of $X$ are disjoint, the Additivity Property
implies
$$
I(f) = i(X, f, S^1_1 - \bar U)  + \sum_{j = 2}^k i(X, f, S^1_j -
\bar U).\leqno(4.2)
$$
The Additivity Property also implies that
$$
I(f_j) = i(S^1_j, f_j, S^1_j - \bar U) + i(S^1_j, f_j,
S^1_j \cap U).\leqno(4.3)
$$
There is a
neighborhood $W_j$ of $(Fix\, f) \cap S^1_j$ in $S^1_j$ such that
$f(\overline W_j) \subseteq S^1_j$.
Thus $f_j(x) = f(x)$ for $x \in W_j$ and therefore, by the
Excision Property,
$$
i(S^1_j, f_j, S^1_j - \overline U) = i(S^1_j, f_j, W_j) = i(X, f,
W_j) = i(X, f, S^1_j - \overline U).\leqno(4.4)
$$

Since $f(\overline U) \subseteq S^1_1$, then $f_1(x) = f(x)$ for all
$x \in \overline U \cap S^1_1$.
There are no fixed points of $f$ in $\overline U$, so
$i(S^1_1, f_1, S^1_1 \cap U) = 0$ and thus $I(f_1) = i(X, f,
S^1_1 - \overline U)$ by (4.3) and (4.4).

For $j \ge 2$, the fact that $f_j(U) = *$ gives us
$i(S^1_j, f_j, S^1_j \cap
U) = 1$ so $I(f_j) = i(X, f, S^1_j - \overline
U) + 1$ by (4.3) and (4.4).
Since $f_j \colon S^1_j \to S^1_j$, the $k = 1$ case
of the argument tells us that $I(f_j) = -deg(f_j) + 1$ for $j =
1, 2, \dots k$.  In particular, $i(X, f, S^1_1 - \overline U) =
-deg(f_1) + 1$
whereas, for $j \ge
2$, we have $i(X, f, S^1_j - \overline U) = -deg(f_j)$.  Therefore, by
(4.2),
$$
I(f) = i(X, f, S^1_1 - \overline U)  + \sum_{j = 2}^k i(X, f, S^1_j -
\overline U) = -\sum^k_{j = 1} deg(f_j) + 1.
$$
This completes the proof of Theorem 1.3.
\qed

\enddocument